\documentclass[11pt]{article}
\usepackage[margin=1in]{geometry}
\usepackage{graphicx,xcolor,cite,mathtools}
\usepackage{amssymb,amsmath,amsthm,mathrsfs,paralist,bm,esint,setspace}
\RequirePackage[colorlinks,citecolor=blue,urlcolor=blue]{hyperref}

\newcommand{\R}{{\mathbb{R}}}
\newcommand{\Z}{{\mathbb{Z}}}
\newcommand{\ZZ}{{\mathrm{Z}}}

\newcommand{\E}{\mathrm{E}}
\newcommand{\N}{\mathrm{N}}

\newcommand{\HH}{\mathfrak{H}}
\renewcommand{\P}{\mathrm{P}}

\newcommand{\e}{\mathrm{e}}

\newcommand{\Var}{\text{\rm Var}}

\input cyracc.def

\title{
Gaussian fluctuation for  Gaussian Wishart matrices of overall correlation\thanks{
	Research supported in part by   FNR grant APOGee (R-AGR-3585-10) at University of Luxembourg.}}
\author{Ivan Nourdin and  Fei Pu\\University of Luxembourg
	}
\date{\today}

\begin{document}
\newtheorem{stat}{Statement}[section]
\newtheorem{proposition}[stat]{Proposition}
\newtheorem*{prop}{Proposition}
\newtheorem{corollary}[stat]{Corollary}
\newtheorem{theorem}[stat]{Theorem}
\newtheorem{lemma}[stat]{Lemma}
\theoremstyle{definition}
\newtheorem{definition}[stat]{Definition}
\newtheorem*{cremark}{Remark}
\newtheorem{remark}[stat]{Remark}
\newtheorem*{OP}{Open Problem}
\newtheorem{example}[stat]{Example}
\newtheorem{nota}[stat]{Notation}

\numberwithin{equation}{section}
\maketitle

\begin{abstract}
In this note, we study  the Gaussian {fluctuations} for the Wishart matrices  $d^{-1}\mathcal{X}_{n, d}\mathcal{X}^{T}_{n, d}$, where 
$\mathcal{X}_{n, d}$ is a $n\times d$ random matrix whose entries are jointly Gaussian and correlated with row and column covariance functions given by $r$ and $s$ respectively such that $r(0)=s(0)=1$. Under the assumptions $s\in \ell^{4/3}(\Z)$ and $\|r\|_{\ell^1(\Z)}< \sqrt{6}/2$, 
we establish the
$\sqrt{n^3/d}$ convergence rate for the Wasserstein distance between a normalization of $d^{-1}\mathcal{X}_{n, d}\mathcal{X}^{T}_{n, d}$ and the corresponding Gaussian ensemble. 
This rate is the same as the optimal one computed in \cite{JL15,BG16,BDER16} for the total variation distance, in the particular case where the Gaussian entries of $\mathcal{X}_{n, d}$ are independent.
Similarly, we obtain the 
$\sqrt{n^{2p-1}/d}$ convergence rate  for the {Wasserstein distance} in the setting of random $p$-tensors of overall correlation.
Our analysis is based on the Malliavin-Stein approach.
\end{abstract}

\bigskip

 \smallskip
\noindent{\it Abbreviated title}: Gaussian fluctuation for  Wishart matrices

\section{Introduction and main result}

Let $\HH$ be a real separable Hilbert space equipped with the inner product $\langle \cdot\, , \cdot\rangle_{\HH}$
 and the Hilbert norm $\left\|\cdot\right\|_{\HH}$, and let $\{e_{ij}: i, j\geq 1\} \subset \HH$ be a family such that 
 \begin{align}\label{covariance}
\langle e_{ij}\,, e_{i'j'}\rangle_{\HH} = r(i-i') \,s(j-j'),
\end{align}
where $s, r: \Z\to \R$ stand for some covariance functions satisfying  $s(0)=r(0)=1$. In particular,  {observe} that $\left\|e_{ij}\right\|_{\HH}=1$ for all $i, j\geq 1$.

Consider the corresponding Gaussian sequence $X_{ij} = X(e_{ij}) \sim \N(0, 1)$ where $X=\{X(h), h\in \HH\}$ is  {an isonormal Gaussian process over 
$\HH$, that is, a centered Gaussian process indexed by $\HH$ such that $\E[X(g)X(h)] = \langle g\,, h\rangle_{\HH}$ for all $g, h\in \HH$}. Let $\mathcal{X}_{n, d}$ be the $n\times d$ random matrix given by 
\begin{align}\label{rm}
\mathcal{X}_{n,d}=(X_{ij})_{1\leq i\leq n, 1\leq j\leq d}=\left(\begin{array}{cccc}X_{11} & X_{12} & \ldots & X_{1d}\\X_{21} & X_{22} & \ldots& X_{2d} \\ \vdots & \vdots & \vdots& \vdots\\ X_{n1} & X_{n2} & \cdots& X_{nd}\end{array}\right).
\end{align}

Our goal is to study the high-dimensional  {fluctuations} of Gaussian Wishart matrices $d^{-1}\mathcal{X}_{n, d}\mathcal{X}^{T}_{n, d}$ by considering a normalized 
version given by
\begin{align}\label{Wishart}
\widetilde{\mathcal{W}}_{n, d}=\left( \widetilde{W}_{ij}\right)_{1\leq i, j\leq n},
\end{align}
where 
\begin{align}\label{Wishart-entry}
\widetilde{W}_{ij} = \frac{1}{\sqrt{d}}\sum_{k=1}^d \left(X_{ik} X_{jk} - r(i-j)\right).
\end{align}
{Since the $e_{ij}$'s are not supposed to be orthogonal, see (\ref{covariance}), it is important to note that the Gaussian entries $X_{ij}$ of $\mathcal{X}_{n, d}$  are fully correlated in general.}

Let $\mathcal{G}_{n, d}^{r, s}= \left(G_{ij}\right)_{1\leq i, j\leq n}$ be the $n\times n$ random symmetric matrix such that the associated random vector $\left(G_{11}, \ldots, G_{1n}, G_{21}, \ldots, G_{2n}, \ldots, G_{n1}, \ldots, G_{nn} \right)$ is Gaussian with mean $0$ and has the same covariance matrix as 
\begin{align*}
\left(\widetilde{W}_{11}, \ldots, \widetilde{W}_{1n}, \widetilde{W}_{21}, \ldots, \widetilde{W}_{2n}, \ldots, \widetilde{W}_{n1}, \ldots, \widetilde{W}_{nn}\right).
\end{align*}

Recall the definition of {\em Wasserstein distance} between two random variables with values in $\mathcal{M}_n(\R)$ (the space of $n\times n$
real matrices): for  {$\mathcal{X}, \mathcal{Y}:\Omega\to \mathcal{M}_n(\R)$ such that $\E\|\mathcal{X}\|_{\rm HS}+\E\|\mathcal{Y}\|_{\rm HS}<\infty$}, 
\begin{align}\label{wd}
d_{\rm Wass} (\mathcal{X}\,, \mathcal{Y}):= \sup\left\{ \E[g(\mathcal{X})] - \E[g(\mathcal{Y})]: \left\|g\right\|_{\rm Lip} \leq 1  \right\},
\end{align}
with
 $$
  \left\|g\right\|_{\rm Lip} := \sup_{\substack{A,B\in \mathcal{M}_n(\R)\\A\neq B}}\frac{|g(A)- g(B)|}{\|A-B\|_{\rm HS}}\quad  {\mbox{for $g:\mathcal{M}_n(\R)\to\R$}},
$$
 {and} $\|\cdot\|_{\rm HS}$  the Hilbert-Schmidt norm on $\mathcal{M}_n(\R)$.

The main result of this paper is the following. 
\begin{theorem}\label{maintheorem}
          Assume
                \begin{align}\label{r}
          \|r\|_{\ell^1(\Z)} < \sqrt{6}/2.
          \end{align}
          Then for all $n, d\geq 1$,
 {
                   \begin{equation} \label{n3/d}
         d_{\rm Wass}\left(\widetilde{\mathcal{W}}_{n,d}\,, \mathcal{G}^{r, s}_{n, d}\right)
\leq\frac{\|r\|^{3/2}_{\ell^1(\Z)}}{3- 2\|r\|^2_{\ell^1(\Z)}} 
            \sqrt{ \frac{32n^3}{d} 
                              \left(\sum_{|k|\leq d}|s(k)|^{4/3}\right)^{3}}.
                                       \end{equation}  
}    
\end{theorem}
\begin{remark}\label{remark1}

\begin{itemize}
\item [(1)]  {We will actually show that 
$$d_{\rm Wass}\left(\widetilde{\mathcal{W}}_{n,d}\,, \mathcal{G}^{r, s}_{n, d}\right)
\leq  
         \frac{\|r\|^{3/2}_{\ell^1(\Z)}}{3- 2\|r\|^2_{\ell^1(\Z)}} 
          \sqrt{\frac{32d}{\sum_{k,\ell=1}^{d}s(k-\ell)^2} \times \frac{n^3}{d} 
                            \left(\sum_{|k|\leq d}|s(k)|^{4/3}\right)^{3}}.$$
This implies (\ref{n3/d}) since $\sum_{k,\ell=1}^{d}s(k-\ell)^2 \geq d\, s(0)^2=d$.
}

\item [(2)]
Under the condition \eqref{r} and if we assume $s\in \ell^{4/3}(\Z)$, then  \eqref{n3/d} leads to 
$$d_{\rm Wass}\left(\widetilde{\mathcal{W}}_{n,d}\,, \mathcal{G}^{r, s}_{n, d}\right) = O(\sqrt{n^3/d}).$$
Hence in this case of overall correlation, $\mathcal{W}_{n,d}$ continues to be close to the Gaussian random matrix $ \mathcal{G}^{r, s}_{n, d}$ as long as $n^3/d \to 0$, 
exactly like in the row independence case in \cite[Theorem 1.2]{NourdinZheng2018}; see also the full independence case considered in 
\cite{BDER16, BG16, JL15}.

\item [(3)] An explicit example of covariance function satisfying  \eqref{r} is 
$r(k)= \e^{-\lambda |k|^{\alpha}}$ for $k\in \Z$, with $1\leq \alpha \leq 2$ and where $\lambda>0$ 
is chosen large enough.

\end{itemize}
\end{remark}

The 
 $\sqrt{n^3/d}$ convergence rate obtained in Theorem \ref{maintheorem} relies on Malliavin calculus and  {Stein's} method, precisely, Proposition \ref{malliavin-stein} below, which  has already been employed in \cite{NourdinZheng2018} to investigate the Gaussian approximation for Wishart matrix in the row independence case ({that is,} $r(k)=1_{k=0}$). In the case of overall correlation, 
it is not  clear if the covariance matrix $C$ in Proposition \ref{malliavin-stein} is invertible or not. {To bypass this problem, the authors of \cite{NourdinZheng2018} made use of the bounds from \cite[Theorem 6.1.2]{NP} and \cite[Theorem 9.3]{NZ17} with, as a price to pay, the necessity to consider} a smoother distance instead of Wasserstein distance; see \cite[Proposition 4.1 and Theorem 4.3]{NourdinZheng2018}.

Fortunately, in the case of overall correlation, we discover that condition \eqref{r} guarantees that the  covariance matrix $C$ in Proposition \ref{malliavin-stein} is  strictly diagonally dominant and hence invertible. Moreover, we can bound the operator norms $\|C\|_{\rm op}$ and $\|C^{-1}\|_{\rm op}$ in terms of $\|r\|_{\ell^1(\Z)}$. Therefore, we are able to apply the inequality \eqref{msinequality} in Proposition \ref{malliavin-stein} to derive the estimate in  \eqref{n3/d}; see the proof of Theorem \ref{maintheorem} in Section \ref{proof1}.

The Malliavin-Stein approach can also be applied to study Gaussian approximation of Wishart $p$-tensors in the case of overall correlation.  In Theorem \ref{theorem:tensor}, we propose a condition on $\|r\|_{\ell^1(\Z)}$ (see \eqref{condition:r} below) under which 
the covariance matrix of the $p$-tensors is invertible. Hence we appeal to the Proposition \ref{malliavin-stein} again and establish the $\sqrt{n^{2p-1}/d}$ convergence rate for the Wasserstein distance between the $p$-tensors; the same as  full independence case considered in
 \cite[Theorem 4.3]{NourdinZheng2018}.
We refer to \cite{BD20, BDT20, M20} for some other recent applications of Malliavin calculus and Stein method in the study of high-dimensional regime of Wishart matrices/tensors.

\section{Preliminaries}

In this section, we collect some elements of Malliavin calculus and  {Stein's} method and refer to \cite{NP} (see also \cite{Nualart, NN}) for more details. 
Recall the isonormal Gaussian process $X=\{X(h), h\in \HH\}$ over a real separable Hilbert space $\HH$ defined on some probability space $(\Omega, \mathcal{F}, \P)$.

For every $p\geq 1$, we let $\mathcal{H}_p$ denote the $p$th Wiener chaos of $X$, that is, the closed linear subspace of $L^2(\Omega)$
generated by the random variables of the form $\{H_p(X(h)), h\in \HH, \|h\|_{\HH}=1\}$, where $H_p$ stands for the $p$th Hermite polynomial. 
The relation that $I_p(h^{\otimes p}) = H_p(X(h))$ for unit vector $h\in \HH$ can be extended to a linear isometry between the symmetric $p$th
tensor product $\HH^{\odot p}$ and the $p$th Wiener chaos $\mathcal{H}_p$.

Consider $f\in \HH^{\odot p}$ and $g\in \HH^{\odot q}$ with $p, q\geq 1$. For $j\in \{0, \ldots, p\wedge q\}$, $f\otimes_{j} g$ denotes the $j$-contraction of $f$ and $g$ (see \cite[Section B.4]{NP} for the precise definition) and  $f\widetilde{\otimes}_{j} g$ stands for the symmetrization of  $f\otimes_{j} g$. For $f\in \HH^{\odot p}$, the Malliavin derivative of $I_p(f)$ is the random element of $\HH$ given by 
$DI_p(f) = pI_{p-1}(f)$ (see \cite[Proposition 2.7.4]{NP}) and we have for $f, g \in  \HH^{\odot p}$,
\begin{align}\label{expectation}
\E\left[p^{-1} \langle DI_p(f)\,, DI_p(g)\rangle_{\HH}\right]= \E[I_p(f)I_p(g)] = p!\langle f\,, g\rangle_{\HH^{\otimes p}}.
\end{align}
Moreover, according to the formula \cite[(6.2.3)]{NP}, for $f, g \in \HH^{\odot p}$,
\begin{align}\label{var}
\Var\left(p^{-1} \langle DI_p(f)\,, DI_p(g) \rangle_{\HH}\right) =p^2 \sum_{j=1}^{p-1}(j-1)!^2{p-1\choose j-1}^4(2p-2j)
\| {f\widetilde{\otimes}_{j} g}\|^2_{\HH^{\otimes (2p-2j)}}.
\end{align}

The following result, the  {so-called} Malliavin-Stein approach, provides a powerful machinery to investigate the normal approximation for the Gaussian Wishart  matrix of overall correlation. 
\begin{proposition}[see {\cite[Corollary 3.6]{NPR10}}]\label{malliavin-stein}
          Fix integers $m\geq 2$ and $1\leq p_1\leq \ldots\leq p_m$. Consider a random 
          vector $F= (F_1, \ldots, F_m) = (I_{p_1}(f_1), \ldots, I_{p_m}(f_m))$
          with $f_j\in \HH^{\odot p_j}$ for each $j$.  On the other hand, let $C$ be an invertible covariance matrix and let $\ZZ\in \N_m(0, C)$. Then
          \begin{align}
          d_{\rm Wass}(F, \ZZ) \leq \|C^{-1}\|_{\rm op}\|C\|^{1/2}_{\rm op}\left(\sum_{1\leq i, j\leq m}\E\left[\left(C_{ij}-p_j^{-1}
          \langle DF_i\,, DF_j\rangle_{\HH}
          \right)^2\right]\right)^{1/2}, \label{msinequality}
          \end{align}
          where $\|\cdot\|_{\rm op}$ denotes the usual operator norm. 
\end{proposition}
 Note that the {\em Wasserstein distance} $d_{\rm Wass}(F, \ZZ) $ between two general $m$-dimensional random \underline{vectors} $F$ and $\ZZ$ is defined as
\begin{align}\label{wd-vector}
d_{\rm Wass}(F\,, \ZZ):=\sup \left\{ \E\left[g(F)\right]  -  \E\left[g(\ZZ)\right] : \|g\|_{\rm Lip}\leq 1\right\},
\end{align}
where $\|g\|_{\rm Lip}$ denotes the usual Lipschitz  constant of a function $g: \R^m\to \R$ {with respect to the Euclidean norm}.

Lemma 2.2 of \cite{NourdinZheng2018} has provided a trick to pass the high-dimensional regime for the full-size  {symmetric} matrix to that of half-matrix. 
Recall that the half matrix $\mathcal{Z}^{\rm half}$ of a $n\times n$ random symmetric matrix $\mathcal{Z}=(Z_{ij})_{1\leq i, j \leq n}$ is the $n(n+1)/2$- dimensional random vector formed by the upper-triangular entries, namely:
\begin{align}\label{half}
\mathcal{Z}^{\rm half} = (Z_{11}, Z_{12}, \ldots, Z_{1n}, Z_{22}, \ldots, Z_{23}, \ldots, Z_{2n}, \ldots, Z_{nn}).
\end{align}
According to  \cite[Lemma 2.2]{NourdinZheng2018}, for two symmetric random matrices $\mathcal{X}, \mathcal{Y} {:\Omega\to} \mathcal{M}_n(\R)$, 
\begin{align}\label{full-half}
d_{\rm Wass}(\mathcal{X}\,, \mathcal{Y}) \leq \sqrt{2}\, d_{\rm Wass}(\mathcal{X}^{\rm half}\,, \mathcal{Y}^{\rm half}),
\end{align}
where the left-hand side Wasserstein distance is defined in \eqref{wd} while the right-hand defined in \eqref{wd-vector}.

Finally, in order to apply Proposition \ref{malliavin-stein} to obtain the rate for the Wasserstein distance stated in Theorem \ref{maintheorem},  we use the product formula for the multiple Wiener-It\^o integrals (see \cite[Theorem 2.7.10]{NP}) to realize the  $(i, j)$th entry of 
$\widetilde{\mathcal{W}}_{n, d}$ defined in \eqref{Wishart-entry} as an element in the second Wiener chaos $\mathcal{H}_2$, namely: 
\begin{align}\label{Wishart-entry-int}
\widetilde{W}_{ij} &= \frac{1}{\sqrt{d}}\sum_{k=1}^d \left(X_{ik} X_{jk} - r(i-j)\right)= \frac{1}{\sqrt{d}}\sum_{k=1}^d \left(I_1(e_{ik}) I_1(e_{jk}) - \langle e_{ik}\,, e_{jk}\rangle_{\HH}\right)\nonumber\\
&= I_2(f_{ij}^{(d)}), 
\end{align}
where 
\begin{align}\label{f_ij}
f_{ij}^{(d)} = \frac{1}{\sqrt{d}}\sum_{k=1}^d e_{ik}\widetilde{\otimes}e_{jk}= \frac{1}{2\sqrt{d}}\sum_{k=1}^d ( e_{ik}\otimes e_{jk} +e_{jk}\otimes e_{ik}).
\end{align}

\section{Proof of Theorem \ref{maintheorem}}\label{proof1}

We prove Theorem \ref{maintheorem} in this section and assume that \eqref{r} holds throughout this section.
We begin to establish the following supporting lemmas.

\begin{lemma}\label{bound}
          For all $n\geq 1$ and for all fixed $(i, j)$ with $1\leq i\leq j \leq n$, 
          \begin{align}\label{e01}
         \sum_{\substack{1\leq u\leq v \leq n \\ (u, v)\neq (i, j)}} \left|r(i-u)r(v-j) + r(i-v)r(u-j)\right|\leq 2 \|r\|^2_{\ell^1(\Z)} -2 < 1. 
          \end{align}
\end{lemma}
\begin{proof}
           The second inequality in \eqref{e01} is clearly true by \eqref{r}. In order to prove the first one, 
           we write 
          \begin{align}\label{sums}
          &\sum_{\substack{1\leq u\leq v \leq n \\ (u, v)\neq (i, j)}} \left|r(i-u)r(v-j) + r(i-v)r(u-j)\right|  \nonumber\\
          &\qquad \qquad \qquad\qquad \qquad \qquad
          \leq \sum_{\substack{1\leq u\leq v \leq n \\ (u, v)\neq (i, j)}} \left|r(i-u)r(v-j) \right| 
          + \sum_{\substack{1\leq u\leq v \leq n \\ (u, v)\neq (i, j)}} \left| r(i-v)r(u-j)\right|.
          \end{align}
          For the first sum on the right-hand side of \eqref{sums}, we have
          \begin{align}\label{e0}
          \sum_{\substack{1\leq u\leq v \leq n \\ (u, v)\neq (i, j)}} \left|r(i-u)r(v-j)\right| &\leq r(0)\sum_{v\neq j}|r(v-j)|
          + \sum_{u\neq i}\sum_{v\in \Z}|r(i-u)|\, |r(v-j)| \nonumber\\
          & = \left(r(0) + \|r\|_{\ell^1(\Z)}\right)\|r\|_{\ell^1(\Z\setminus \{0\})} \nonumber\\
          &= \|r\|^2_{\ell^1(\Z)}-1 .
          \end{align}
          For the second sum on the right-hand side of \eqref{sums},
          \begin{align}\label{e1}
            \sum_{\substack{1\leq u\leq v \leq n \\ (u, v)\neq (i, j)}} \left| r(i-v)r(u-j)\right| \leq |r(i-j)|\sum_{v\neq j}|r(i-v)| 
            +             \sum_{\substack{1\leq u\leq v \leq n \\ u\neq i}} \left| r(i-v)r(u-j)\right|
          \end{align}
          We observe that 
          \begin{align}\label{e2}
           |r(i-j)|\sum_{v\neq j}|r(i-v)| \leq 
           \begin{cases}
           r(0)\|r\|_{\ell^{1}(\Z\setminus \{0\})}=  \|r\|_{\ell^1(\Z)}-1, & \text{if}\,\, i=j, \\
           \frac{1}{2}\|r\|_{\ell^{1}(\Z\setminus \{0\})}           \|r\|_{\ell^{1}(\Z)} =   \|r\|_{\ell^1(\Z)}(\|r\|_{\ell^1(\Z)}-1)/2, & \text{if}\,\, i<j.
           \end{cases}
          \end{align}
          Moreover,  since $i\leq j$, 
                    \begin{align}\label{e3}
           \sum_{\substack{1\leq u\leq v \leq n \\ u\neq i}} \left| r(i-v)r(u-j)\right| &=          
             \sum_{\substack{1\leq u\leq v \leq n \\ u< i}} \left| r(i-v)r(u-j)\right| 
             +            \sum_{\substack{i<u\leq v \leq n }} \left| r(i-v)r(u-j)\right| \nonumber \\
             &\leq  \|r\|_{\ell^1(\Z)} \sum_{1\leq u<i}|r(u-j)| + \sum_{i< v\leq n}|r(i-v)| \|r\|_{\ell^1(\Z)} \nonumber \\
             &\leq  \|r\|_{\ell^1(\Z)}\|r\|_{\ell^1(\Z\setminus\{0\})} =\|r\|_{\ell^1(\Z)}(\|r\|_{\ell^1(\Z)}-1).
          \end{align}
          By \eqref{e1}--\eqref{e3}, it yields that 
          \begin{align}\label{e4}
            \sum_{\substack{1\leq u\leq v \leq n \\ (u, v)\neq (i, j)}} \left| r(i-v)r(u-j)\right|
             \leq (\|r\|_{\ell^1(\Z)}-1)\left(\left(1\vee (\|r\|_{\ell^1{(\Z)}}/2)\right) +\|r\|_{\ell^1{(\Z)}} \right).
                      \end{align}
           Therefore, we combine \eqref{e0} and \eqref{e4} to obtain
                     \begin{align}\label{e011}
         & \sum_{\substack{1\leq u\leq v \leq n \\ (u, v)\neq (i, j)}} \left|r(i-u)r(v-j) + r(i-v)r(u-j)\right| \nonumber\\
         &\qquad \qquad \qquad \qquad \qquad\leq  (\|r\|_{\ell^1(\Z)}-1)\left(\left(1\vee (\|r\|_{\ell^1{(\Z)}}/2)\right) +2\|r\|_{\ell^1{(\Z)}} +1\right),
          \end{align}
          Elementary calculation on quadratic inequality shows that
          \begin{align*}
           (\|r\|_{\ell^1(\Z)}-1)\left(\left(1\vee (\|r\|_{\ell^1{(\Z)}}/2)\right) +2\|r\|_{\ell^1{(\Z)}} +1\right)< 1\Leftrightarrow \|r\|_{\ell^1(\Z)} < \sqrt{6}/2,
          \end{align*}
          whence 
           \begin{align}\label{iff}
           (\|r\|_{\ell^1(\Z)}-1)\left(\left(1\vee (\|r\|_{\ell^1{(\Z)}}/2)\right) +2\|r\|_{\ell^1{(\Z)}} +1\right)=2 \|r\|^2_{\ell^1(\Z)} -2
          \end{align}
           provided $\|r\|_{\ell^1(\Z)} <  \sqrt{6}/2$.
           
           Therefore, under the assumption \eqref{r},  the estimate \eqref{e01} follows from \eqref{e011} and \eqref{iff}.
\end{proof}

Recall the $n\times n$ symmetric Gaussian random matrix  $\mathcal{G}_{n, d}^{r, s}= \left(G_{ij}\right)_{1\leq i, j\leq n}$  {from} \eqref{Wishart-entry}. 
Let $C$ denote the covariance matrix of $(\mathcal{G}^{r, s}_{n, d})^{\text{half}}$.  Notice that the matrix norms 
$\|C\|_1$ and $\|C\|_{\infty}$ of the symmetric matrix $C$ are equal and given by
\begin{align}\label{1infty}
\|C\|_1=\|C\|_{\infty}=\sup_{1\leq i\leq j\leq n} \sum_{\substack{1\leq u\leq v\leq n}} \left|\E[G_{ij}G_{uv}]\right|.
\end{align}

\begin{lemma}\label{invertible}
          The matrix $C$ is invertible and the following estimates on operator norms hold:
          \begin{align}\label{operatornorm}
          \|C^{-1}\|_{\rm op} \leq            \frac{d}{\sum_{k,\ell=1}^{d}s(k-\ell)^2} \left(3- 2\|r\|^2_{\ell^1(\Z)}\right)^{-1}  
          \quad \text{and}\quad \|C\|_{\rm op} \leq   \frac{2 \|r\|^2_{\ell^1(\Z)} }{d}\sum_{k,\ell=1}^{d}s(k-\ell)^2.
          \end{align}
\end{lemma}
\begin{proof}
          According to  \cite[(4.3)]{NourdinZheng2018}, the entries of $C$ are given by 
          \begin{align}\label{cov}
         \E[G_{ij}G_{uv}] =\frac{ \left(r(i-u)r(v-j) + r(i-v)r(u-j)\right) }{d}\sum_{k,\ell=1}^{d}s(k-\ell)^2
         \end{align}
         for $1\leq i \leq j\leq n$ and $1\leq u\leq v\leq n$. 
         Letting $(u, v)= (i, j)$ in \eqref{cov}, we obtain the following lower and upper bounds on 
         the diagonal entries of $C$:
         \begin{align}\label{supinf}
          \frac{1}{d}\sum_{k,\ell=1}^{d}s(k-\ell)^2 = \inf_{\substack{1\leq i\leq j\leq n}}\E[G_{ij}^2] \leq  \sup_{\substack{1\leq i\leq j\leq n}}\E[G_{ij}^2]
          = \frac{2}{d}\sum_{k,\ell=1}^{d}s(k-\ell)^2.
         \end{align}
         Moreover, under the condition \eqref{r}, we have 
          \begin{align}\label{1infty-estimate}
         \|C\|_1&=\|C\|_{\infty}=\sup_{1\leq i\leq j\leq n} \sum_{\substack{1\leq u\leq v\leq n}} \left|\E[G_{ij}G_{uv}]\right| 
         \leq  \frac{2 \|r\|^2_{\ell^1(\Z)} }{d}\sum_{k,\ell=1}^{d}s(k-\ell)^2 
         \end{align}
         thanks to \eqref{supinf} and Lemma \ref{bound}.  Hence the second inequality in \eqref{operatornorm} follows from \eqref{1infty-estimate} and 
         the H\"{o}lder's inequality 
         for matrix norms: $\|C\|_{\text{op}}\leq \sqrt{\|C\|_1\|C\|_{\infty}}$ (see \cite[Theorem 4.3.1]{Ser10}).
         
         {
         Furthermore,  by \eqref{cov} and \eqref{supinf} in the first equality,
         \begin{eqnarray}\label{SDD}
       &&  \inf_{\substack{1\leq i\leq j\leq n}}\left(\E[G_{ij}^2] - \sum_{\substack{1\leq u\leq v\leq n \\ (u, v)\neq (i, j)}} 
         \left|\E[G_{ij}G_{uv}]\right|\right)  \geq   \inf_{\substack{1\leq i\leq j\leq n}}\E[G_{ij}^2]   -  
 \sup_{\substack{1\leq i\leq j\leq n}} \sum_{\substack{1\leq u\leq v\leq n \\ (u, v)\neq (i, j)}} 
         \left|\E[G_{ij}G_{uv}]\right|   \notag\\
&  =  &
        \frac{1}{d}\sum_{k,\ell=1}^{d}s(k-\ell)^2 \left(
1- \sup_{\substack{1\leq i\leq j\leq n}} \sum_{\substack{1\leq u\leq v \leq n \\ (u, v)\neq (i, j)}} \left|r(i-u)r(v-j) + r(i-v)r(u-j)\right|
\right) \notag\\
 &  \geq   &
        \frac{1}{d}\sum_{k,\ell=1}^{d}s(k-\ell)^2 \left(3- 2\|r\|^2_{\ell^1(\Z)}\right)>0,\,\,\mbox{ by Lemma \ref{bound} and under  \eqref{r}},
        \end{eqnarray}
}
        which implies that the symmetric matrix $C$ is strictly diagonally dominant and hence invertible.
        
        Now we apply \cite[Corollary 2]{Varah1975} and \eqref{SDD} to see that 
        \begin{align}\label{operator-lower}
        \|C^{-1}\|_{\rm op}^{-1}& \geq \inf_{\substack{1\leq i\leq j\leq n}}\left(\E[G_{ij}^2] - \sum_{\substack{1\leq u\leq v\leq n \\ (u, v)\neq (i, j)}} 
         \left|\E[G_{ij}G_{uv}]\right|\right)   \nonumber\\
         &              \geq           \frac{1}{d}\sum_{k,\ell=1}^{d}s(k-\ell)^2 \left(3- 2\|r\|^2_{\ell^1(\Z)}\right),
        \end{align}
        which proves the first inequality in \eqref{operatornorm}.
\end{proof}

We are now ready to prove Theorem \ref{maintheorem}.

\begin{proof}[Proof of Theorem \ref{maintheorem}]
          Recall that $(\widetilde{\mathcal{W}}_{n,d})^{\text{half}}$ is the half matrix of the Wishart matrix $\widetilde{\mathcal{W}}_{n,d}$
          defined in \eqref{Wishart} and \eqref{Wishart-entry}, whose entries can be represented as the elements in the second Wiener chaos $
          \mathcal{H}_2$; see \eqref{Wishart-entry-int} and \eqref{f_ij}. By Lemma \ref{invertible}, 
          the covariance matrix of $(\widetilde{\mathcal{W}}_{n,d})^{\text{half}}$ is invertible and hence we can apply Proposition \ref{malliavin-stein} with
          $m=n(n+1)/2$, $p_1=\ldots =p_m=2$ and $F=(\widetilde{\mathcal{W}}_{n,d})^{\text{half}}$.
          Indeed, we have 
          \begin{align}
         d_{\text{Wass}}\left((\widetilde{\mathcal{W}}_{n,d})^{\text{half}}\,, 
         (\mathcal{G}^{r, s}_{n, d})^{\text{half}}\right)&= d_{\text{Wass}}\left(F\,, 
         (\mathcal{G}^{r, s}_{n, d})^{\text{half}}\right) \nonumber\\
        & \leq \|C^{-1}\|_{\rm op}\|C\|^{1/2}_{\rm op}\left(\sum_{1\leq i, j\leq m}\E\left[\left(C_{ij}-\frac12
          \langle DF_i\,, DF_j\rangle_{\HH}
          \right)^2\right]\right)^{1/2}.\nonumber
         \end{align}
         Using the identities \eqref{expectation} and \eqref{var}, the proceeding yields that
         \begin{align}\label{wass-estimate}
         d_{\text{Wass}}\left((\widetilde{\mathcal{W}}_{n,d})^{\text{half}}\,, 
         (\mathcal{G}^{r, s}_{n, d})^{\text{half}}\right)  & \leq \|C^{-1}\|_{\rm op}\|C\|^{1/2}_{\rm op}         \left(\sum_{\substack{1\leq i\leq j\leq n\\ 1\leq p\leq q\leq n}}\Var\left(\frac{1}{2}
         \langle D\widetilde{W}_{ij}\,,   D\widetilde{W}_{pq}\rangle_{\HH} \right)
         \right)^{1/2}\nonumber\\
         &= \|C^{-1}\|_{\rm op}\|C\|^{1/2}_{\rm op}
                  \left(8\sum_{\substack{1\leq i\leq j\leq n\\ 1\leq p\leq q\leq n}}
           \|f_{ij}^{(d)} {\widetilde{\otimes}_1} f_{pq}^{(d)}  \|^2_{\HH^{\otimes 2}}       \right)^{1/2} \nonumber\\
           &\leq  \frac{\|r\|_{\ell^1(\Z)}}{3- 2\|r\|^2_{\ell^1(\Z)}} \sqrt{\frac{16d}{\sum_{k,\ell=1}^{d}s(k-\ell)^2}
                             \sum_{\substack{1\leq i\leq j\leq n\\ 1\leq p\leq q\leq n}}
           \|f_{ij}^{(d)} \otimes_1 f_{pq}^{(d)}  \|^2_{\HH^{\otimes 2}}     }, 
         \end{align}
         where the second inequality follows from \eqref{operatornorm} {and the fact that $\|\widetilde{h}\|_{\HH^{\otimes r}}\leq \|h\|_{\HH^{\otimes r}}$}.
         
        It remains to estimate the last term in \eqref{wass-estimate}.  Appealing to \cite[(4.7)]{NourdinZheng2018},
        \begin{align}\label{4.7}
        \|f_{ij}^{(d)} \otimes_1 f_{pq}^{(d)}  \|^2_{\HH^{\otimes 2}} \leq \frac{\mathfrak{X}_{i,j,p,q}}{16d} \left(\sum_{|k|\leq d}|s(k)|^{4/3}\right)^3,
        \end{align}
        where $\mathfrak{X}_{i,j,p,q}$ is a sum of sixteen terms given by the expression below $(4.7)$ in \cite{NourdinZheng2018}. Moreover, we have
        \begin{align}\label{X20}
       \mathfrak{X}_{i,j,p,q}\leq 7|r(j-q)| +5|r(p-j)| + 3 |r(i-q)| + |r(i-p)|,
        \end{align} 
        which together with \cite[(4.12)]{NourdinZheng2018} implies that 
        \begin{align}\label{X2}
       \sum_{\substack{1\leq i\leq j\leq n\\ 1\leq p\leq q\leq n}}\mathfrak{X}_{i,j,p,q}\leq 16n^3\|r\|_{\ell^1(\Z)}.
        \end{align} 
        Taking into account \eqref{wass-estimate}, \eqref{4.7} and \eqref{X2},  we conclude that 
             \begin{align}\label{wass-estimate2}
         d_{\text{Wass}}\left((\widetilde{\mathcal{W}}_{n,d})^{\text{half}}\,, 
         (\mathcal{G}^{r, s}_{n, d})^{\text{half}}\right)  
                    &\leq  \frac{\|r\|^{3/2}_{\ell^1(\Z)}}{3- 2\|r\|^2_{\ell^1(\Z)}} \sqrt{\frac{16d}{\sum_{k,\ell=1}^{d}s(k-\ell)^2}
                            \times \frac{n^3}{d}  \left(\sum_{|k|\leq d}|s(k)|^{4/3}\right)^{3}} .
                                       \end{align}
         
          Finally, we combine \eqref{wass-estimate2} with \eqref{full-half} to obtain the estimate in Remark \ref{remark1}(1), 
          which leads to \eqref{n3/d}
 \end{proof}

\section{Random $p$-tensors} \label{section:tensor}

The result of Theorem \ref{maintheorem} for Gaussian Wishart matrix can be extended to random $p$-tensors ($p\geq 2$). We first introduce some notations of $p$-tensors.  Let $\mathbb{X}_i = (X_{1i}, \ldots, X_{ni})^{T}$ be the $i$th column of the random matrix $\mathcal{X}_{n,d}$ defined in 
\eqref{rm}. We write 
\begin{align*}
\mathbb{X}_i = \sum_{j=1}^nX_{ji} \,\varepsilon_j,
\end{align*}
where $\{\varepsilon_j, j=1, \ldots, n\}$ is the canonical basis of $\R^n$. Then the $p$-tensor product of $\mathbb{X}_i$ is given by
\begin{align*}
\mathbb{X}_i^{\otimes p} = \left(\sum_{j=1}^nX_{ji} \,\varepsilon_j\right)^{\otimes p} = \sum_{j_1, \ldots, j_p=1}^{n} \left(\prod_{k=1}^{p}X_{j_ki}\right)
\varepsilon_{j_1}\otimes \ldots \otimes \varepsilon_{j_p}
\end{align*}
so that
\begin{align*}
\frac{1}{\sqrt{d}}\sum_{i=1}^d\mathbb{X}_i^{\otimes p}  = \sum_{j_1, \ldots, j_p=1}^{n} \frac{1}{\sqrt{d}}\sum_{i=1}^d\left(\prod_{k=1}^{p}X_{j_ki}\right)
\varepsilon_{j_1}\otimes \ldots \otimes \varepsilon_{j_p}.
\end{align*}
A repeated application of the product formula  for the multiple Wiener-It\^o integrals (see \cite[Theorem 2.7.10]{NP}) ensures that 
\begin{align*}
\prod_{k=1}^{p}X_{j_ki} = \prod_{k=1}^{p}I_1(e_{j_ki}) = I_{p}\left(\text{sym}\left(e_{j_1i}\otimes \cdots \otimes e_{j_pi}\right)\right) + \text{lower order terms},
\end{align*}
where sym denotes the canonical symmetrization.

Analogous to the choice of $\widetilde{W}_{ij}$ defined in \eqref{Wishart-entry-int} and \eqref{f_ij}, in the case of overall correlation, we consider the following normalized version of $p$-tensor of $\mathcal{X}_{n, d}$: 
\begin{align*}
\left(\widetilde{\mathbf{Y}}_{\mathbf{j}} =I_p(f_\mathbf{j}^{(d)}),\,  \mathbf{j}= (j_1, \ldots, j_p)\in \{1, \ldots, n\}^p \right),
\end{align*}
where
\begin{align}\label{fd}
f_\mathbf{j}^{(d)}= \frac{1}{\sqrt{d}}\sum_{k=1}^d\text{sym}\left(e_{j_1k}\otimes \cdots \otimes e_{j_pk}\right).
\end{align}
Moreover, we remove the diagonal terms and focus on the Gaussian approximation of 
\begin{align}\label{Y_nd}
\widetilde{\mathcal{Y}}_{n, d} = \left(\widetilde{\mathbf{Y}}_{\mathbf{j}} =I_p(f_\mathbf{j}^{(d)}),\,  \mathbf{j}\in \Delta_p \right),
\end{align}
where $f_\mathbf{j}^{(d)}$ is defined in \eqref{fd} and $\Delta_p = \{(j_1, \ldots, j_p)\in \{1, \ldots, n\}^p: j_1, \ldots, j_p\, \, \text{are mutually distinct}\}$.

The following result extends the Gaussian approximation of random $p$-tensors of full independence in \cite[Theorem 4.6]{NourdinZheng2018} to the case of overall correlation. 
\begin{theorem}\label{theorem:tensor}
           Let $\widetilde{\mathcal{Y}}_{n, d}$  be defined in \eqref{Y_nd} and $\mathbf{Z}= (\mathbf{Z}_{\mathbf{j}}: \mathbf{j}\in \Delta_p)$
             a centered Gaussian vector in $\R^{p!{n\choose p}}$ 
          which has the same covariance matrix as $\widetilde{\mathcal{Y}}_{n, d}$.           Assume that the covariance function $r$ satisfies 
          \begin{align}\label{condition:r}
                        \left(1 -   \left(\|r\|_{\ell^1(\Z)} -1\right)
               (p!  \|r\|^{p-1}_{\ell^1(\Z)} + (p!-1)/2)\right)>0.
          \end{align}  
          Then there exists a positive constant $C_p$ depending on $p$ and $\|r\|_{\ell^1(\Z)}$ (see 
          \eqref{gauss2} below) such that 
          \begin{align}\label{gauss}
          d_{\text{Wass}}\left(\widetilde{\mathcal{Y}}_{n, d}\,, \mathbf{Z}\right) \leq 
          C_p                     \sqrt{
          \frac{d} {\left|\sum_{k,\ell =1}^d s(k-\ell)^p\right|}
          \left(\sum_{|k|\leq d}|s(k)|^{4/3}\right)^{3}\frac{n^{2p-1}}{d}  }.
\end{align}
\end{theorem}
\begin{remark}
          \begin{itemize}
          \item [(1)] Note that the above Wasserstein distance is for $\R^{p!{n\choose p}}$-valued random vectors;   as defined in \eqref{wd-vector}.
          \item [(2)] The estimate in \eqref{gauss} is trivial if $\sum_{k,\ell =1}^d s(k-\ell)^p=0$. Hence we assume
           $\sum_{k,\ell =1}^d s(k-\ell)^p\neq 0$ in the following.
           \item [(3)] Under the condition \eqref{condition:r}, if we assume $s\in \ell^{4/3}(\Z)$
            and in addition that $p$ is even or $s(k)\geq 0$ for all $k\in\Z$,
           \eqref{gauss} leads to $d_{\text{Wass}}(\widetilde{\mathcal{Y}}_{n, d}\,, \mathbf{Z}) = O(\sqrt{n^{2p-1}/d})$; 
           the same as the full independence case considered in \cite[Theorem 4.6]{NourdinZheng2018}.
          \end{itemize}
\end{remark}

Similar to the proof of Theorem \ref{maintheorem}, we reduce the estimate of Wasserstein distance between $\widetilde{\mathcal{Y}}_{n, d}$ and $\mathbf{Z}$ to that of their "half matrices",  given by the following $\R^{{n\choose p}}$-valued random vectors
\begin{align}\label{Yup}
\widetilde{\mathcal{Y}}_{n, d}^{\uparrow}= \left(\widetilde{\mathbf{Y}}_{\mathbf{j}} =I_p(f_\mathbf{j}^{(d)}),\,  \mathbf{j}\in \Delta^{\uparrow}_p \right)
\quad \text{and} \quad \mathbf{Z}^{\uparrow}= \left(\mathbf{Z}_{\mathbf{j}}: \mathbf{j}\in \Delta^{\uparrow}_p\right),
\end{align}
where $\Delta_{p}^{\uparrow} = \{\mathbf{j}\in \{1, \ldots, n\}^p: j_1<j_2< \ldots < j_p\}$.
Denote $S(p)$ the collection of all permutations of $\{1, \ldots, p\}$.

\begin{lemma}\label{invertible-tensor}
          Let $\widetilde{C}$ be the covariance matrix of $\widetilde{\mathcal{Y}}_{n, d}^{\uparrow}$. Under the condition \eqref{condition:r}, $\widetilde{C}$ is 
          invertible and we have 
          \begin{align}\label{op}
          \|\widetilde{C}^{-1}\|_{\rm op} \leq    d \left|\sum_{k,\ell =1}^d s(k-\ell)^p\right|^{-1}  
              \left(1 -   \left(\|r\|_{\ell^1(\Z)} -1\right)
               \left(p!  \|r\|^{p-1}_{\ell^1(\Z)} + (p!-1)/2\right)\right)^{-1}
          \end{align}
          and 
          \begin{align}\label{op2}
          \|\widetilde{C}\|_{\rm op} \leq   \frac{1}{d} \left|\sum_{k,\ell =1}^d s(k-\ell)^p\right|              \left(1 +   \left(\|r\|_{\ell^1(\Z)} -1\right)
               \left(p!  \|r\|^{p-1}_{\ell^1(\Z)} + (p!-1)/2\right)\right).
          \end{align}
\end{lemma}
\begin{proof}
          The proof is similar to that of Lemma \ref{invertible}. We will see that the condition \eqref{condition:r}
           guarantees that the symmetric matrix $\widetilde{C}$ is 
          strictly diagonally dominant and hence invertible. We first compute the entries of $\widetilde{C}$. 
          For $\mathbf{j}= (j_1, \ldots, j_p),\, \mathbf{j}'=(j'_1, \ldots, j'_p) \in \Delta_{p}^{\uparrow} $, using \eqref{Y_nd}, \eqref{fd} and isometry, 
                    \begin{align} \label{entry}
          \E[\widetilde{\mathbf{Y}}_{\mathbf{j}}\widetilde{\mathbf{Y}}_{\mathbf{j}'}]&=  \E[ I_p(f_\mathbf{j}^{(d)}) I_p(f_{\mathbf{j}'}^{(d)})] \nonumber\\
          & =\frac{p!}{d}\sum_{k,\ell =1}^d \left\langle \text{sym}\left(e_{j_1k}\otimes \cdots \otimes e_{j_pk}\right)\,, \text{sym}\left(e_{j'_1\ell}\otimes 
          \cdots    \otimes e_{j'_p\ell}\right)\right\rangle_{\HH^{\otimes p}}\nonumber\\
                    & =\frac{1}{p!d}\sum_{k,\ell =1}^d \sum_{\sigma, \tau\in S(p)}
                    \left\langle e_{j_{\sigma(1)}k}\otimes \cdots \otimes e_{j_{\sigma(p)}k}\,, e_{j'_{\tau(1)}\ell}\otimes \cdots \otimes 
                    e_{j'_{\tau(p)}\ell}\right\rangle_{\HH^{\otimes p}}\nonumber\\
                    & =\frac{1}{p!d}\sum_{k,\ell =1}^d s(k-\ell)^p \sum_{\sigma, \tau\in S(p)}
                    \prod_{m=1}^pr(j_{\sigma(m)}-{j'_{\tau(m)}}).
          \end{align}
          As a consequence of \eqref{entry}, we have the following lower and upper bounds on the diagonal entries of $\widetilde{C}$: 
          for all  $\mathbf{j} \in \Delta_{p}^{\uparrow} $
          \begin{align}\label{dia}
                     \E[\widetilde{\mathbf{Y}}_{\mathbf{j}}^2]&\geq \frac{1}{p!d} \left|\sum_{k,\ell =1}^d s(k-\ell)^p\right|
                     \left(p! r(0)^p - \sum_{\substack{\sigma, \tau\in S(p) \\\sigma\neq \tau}} 
                       \prod_{m=1}^p\left|r(j_{\sigma(m)}-{j_{\tau(m)}})\right|\right)\nonumber\\
                     & \geq \frac{1}{p!d} \left|\sum_{k,\ell =1}^d s(k-\ell)^p\right|
                     \left(p!  -((p!)^2 - p!)\frac{1}{2}\|r\|_{\ell^1(\Z\setminus \{0\})}\right) \nonumber \\
                     &= \frac{1}{d} \left|\sum_{k,\ell =1}^d s(k-\ell)^p\right|
                     \left(1 -(p! - 1)\frac{1}{2}\left(\|r\|_{\ell^1(\Z)}-1\right)\right),
                  \end{align}
          and similarly
                    \begin{align}\label{diaupper}
                     \E[\widetilde{\mathbf{Y}}_{\mathbf{j}}^2]&\leq  \frac{1}{d} \left|\sum_{k,\ell =1}^d s(k-\ell)^p\right|
                     \left(1+ (p! - 1)\frac{1}{2}\left(\|r\|_{\ell^1(\Z)}-1\right)\right).
                  \end{align}

                  Moreover,  the identity \eqref{entry} implies the following estimate on the off-diagonal entries of $\widetilde{C}$:
                   for all  $\mathbf{j} \in \Delta_{p}^{\uparrow} $,
              \begin{align}\label{offdia}
              \sum_{\substack{ \mathbf{j}' \in \Delta_{p}^{\uparrow} \\ \mathbf{j}'\neq \mathbf{j}}}
             \left|          \E[\widetilde{\mathbf{Y}}_{\mathbf{j}}\widetilde{\mathbf{Y}}_{\mathbf{j}'}]\right|& \leq
              \frac{1}{p!d} \left|\sum_{k,\ell =1}^d s(k-\ell)^p\right|
             \sum_{\sigma, \tau\in S(p)}
              \sum_{\substack{ \mathbf{j}' \in \Delta_{p}^{\uparrow} \\ \mathbf{j}'\neq \mathbf{j}}} 
               \prod_{m=1}^p\left|r(j_{\sigma(m)}-{j'_{\tau(m)}})\right|\nonumber \\
               & \leq \frac{1}{p!d} \left|\sum_{k,\ell =1}^d s(k-\ell)^p\right| (p!)^2 \|r\|_{\ell^1(\Z\setminus \{0\})} \|r\|^{p-1}_{\ell^1(\Z)}\nonumber\\
               & = \frac{p!}{d} \left|\sum_{k,\ell =1}^d s(k-\ell)^p\right|  \left(\|r\|_{\ell^1(\Z)} -1\right) \|r\|^{p-1}_{\ell^1(\Z)}.
              \end{align}
              Therefore, we obtain from \eqref{dia} and \eqref{offdia} that
              \begin{align}\label{alphap}
             & \inf_{\mathbf{j} \in \Delta_{p}^{\uparrow}}\left(\E[\widetilde{\mathbf{Y}}_{\mathbf{j}}^2] 
            - \sum_{\substack{ \mathbf{j}' \in \Delta_{p}^{\uparrow}, \,\mathbf{j}'\neq \mathbf{j}}}
           \left|          \E[\widetilde{\mathbf{Y}}_{\mathbf{j}}\widetilde{\mathbf{Y}}_{\mathbf{j}'}]\right|\right) \nonumber\\
            &\qquad \qquad \qquad \geq   \frac{1}{d} \left|\sum_{k,\ell =1}^d s(k-\ell)^p\right|  
              \left(1 -   \left(\|r\|_{\ell^1(\Z)} -1\right)
               \left(p!  \|r\|^{p-1}_{\ell^1(\Z)} + (p!-1)/2\right)\right)>0,
                         \end{align}
            thanks to \eqref{condition:r}. Hence the symmetric matrix $\widetilde{C}$ is 
          strictly diagonally dominant and invertible. One more appeal to \cite[Corollary 2]{Varah1975} yields that 
          \begin{align}\label{op_lower}
          \|\widetilde{C}^{-1}\|_{\rm op}^{-1} \geq  \inf_{\mathbf{j} \in \Delta_{p}^{\uparrow}}\left(\E[\widetilde{\mathbf{Y}}_{\mathbf{j}}^2] 
            - \sum_{\substack{ \mathbf{j}' \in \Delta_{p}^{\uparrow}, \,\mathbf{j}'\neq \mathbf{j}}}
           \left|          \E[\widetilde{\mathbf{Y}}_{\mathbf{j}}\widetilde{\mathbf{Y}}_{\mathbf{j}'}]\right|\right), 
          \end{align}
          which together with \eqref{alphap} proves \eqref{op}. 
          
          Furthermore, we deduce from \eqref{diaupper} and \eqref{offdia} that 
          \begin{align}\label{betap}
           & \|\widetilde{C}\|_{1}=\|\widetilde{C}\|_{\infty}=\sup_{\mathbf{j} \in \Delta_{p}^{\uparrow}}\sum_{\substack{ \mathbf{j}' \in \Delta_{p}^{\uparrow} }}
           \left|          \E[\widetilde{\mathbf{Y}}_{\mathbf{j}}\widetilde{\mathbf{Y}}_{\mathbf{j}'}]\right| \nonumber\\
           &   \qquad\qquad\leq  \frac{1}{d} \left|\sum_{k,\ell =1}^d s(k-\ell)^p\right|              \left(1 +   \left(\|r\|_{\ell^1(\Z)} -1\right)
               \left(p!  \|r\|^{p-1}_{\ell^1(\Z)} + (p!-1)/2\right)\right).
               \end{align}
         Hence \eqref{op2} follows from \eqref{betap} and the H\"{o}lder's inequality 
         for matrix norms.
\end{proof}

\begin{lemma}\label{squarenorm}
          Let $f_{\mathbf{j}}^{(d)}, \mathbf{j}\in \Delta_p^{\uparrow}$ be defined in \eqref{fd}.
          For all  $1\leq q\leq p-1$,
          \begin{align}\label{2p-1d}
           \sum_{\mathbf{j}, \mathbf{j}' \in \Delta_{p}^{\uparrow} }\left\|f_{\mathbf{j}}^{(d)} \otimes_{q}
            f_{\mathbf{j}'}^{(d)}\right\|^2_{\HH^{\otimes 2p-2q}} \leq \|r\|_{\ell^1(\Z)}\frac{n^{2p-1}}{d }\left(\sum_{|k|\leq d}|s(k)|^{4/3}\right)^{3}.
          \end{align}
\end{lemma}
\begin{proof}
          We first compute the norm on the left-hand side of \eqref{2p-1d}. 
          For $\mathbf{j}= (j_1, \ldots, j_p)$ and $\mathbf{j}'= (j'_1, \ldots, j'_p)$,
                    using the definition of $f_{\mathbf{j}}^{(d)}$ and $f_{\mathbf{j}'}^{(d)}$,
          \begin{align}\label{tensor-id}
          f_{\mathbf{j}}^{(d)} \otimes_{q} f_{\mathbf{j}'}^{(d)} &= 
          \frac{1}{d}\sum_{k,\ell=1}^d \text{sym}\left(e_{j_1k}\otimes \cdots \otimes e_{j_pk}\right) \otimes_q 
          \text{sym}\left(e_{j'_1\ell}\otimes \cdots \otimes e_{j'_p\ell}\right)\nonumber\\
          &=\frac{1}{d(p!)^2}\sum_{k,\ell=1}^d \sum_{\sigma, \tau\in S(p)}
          \left(e_{j_{\sigma(1)}k}\otimes \cdots \otimes e_{j_{\sigma(p)}k}\right) \otimes_q
          \left(e_{j'_{\tau(1)}\ell}\otimes \cdots \otimes e_{j'_{\tau(p)}\ell}\right)\nonumber\\
           &=\frac{1}{d(p!)^2}\sum_{k,\ell=1}^d\sum_{\sigma, \tau\in S(p)}
          e_{j_{\sigma(q+1)}k}\otimes \cdots \otimes e_{j_{\sigma(p)}k} \otimes
          e_{j'_{\tau(q+1)}\ell}\otimes \cdots \otimes e_{j'_{\tau(p)}\ell}\nonumber\\
          & \qquad\qquad \qquad \times s(k-\ell)^q \prod_{m=1}^{q}r(j_{\sigma(m)} -j'_{\tau(m)}).
          \end{align}
          Now we take the square norm and it yields that 
          \begin{align}
           &\quad \left\|f_{\mathbf{j}}^{(d)} \otimes_{q} f_{\mathbf{j}'}^{(d)}\right\|^2_{\HH^{\otimes 2p-2q}} \nonumber\\
           &=\frac{1}{d^2(p!)^4}
           \sum_{k,k',\ell, \ell'=1}^d\sum_{\sigma,\sigma', \tau, \tau'\in S(p)}s(k-\ell)^q s(k'-\ell')^q
            \prod_{m=1}^{q}r(j_{\sigma(m)} -j'_{\tau(m)})r(j_{\sigma'(m)} -j'_{\tau'(m)})\nonumber\\
            & \qquad \qquad \quad \times s(k-k')^{p-q}s(\ell-\ell')^{p-q} 
            \prod_{m=q+1}^{p}r(j_{\sigma(m)} - j_{\sigma'(m)})r(j'_{\tau(m)} - j'_{\tau'(m)}). \label{square}
            \end{align}
                      For $1\leq q\leq p-1$, since $|s(k)| \leq 1$ for all $k\in \Z$,
          \begin{align}\label{4/3}
          &\frac{1}{d^2}\sum_{k,k',\ell, \ell'=1}^d \left|s(k-\ell)^q s(k'-\ell')^q s(k-k')^{p-q}s(\ell-\ell')^{p-q} \right|\nonumber\\
          & \quad \leq    \frac{1}{d^2}\sum_{k,k',\ell, \ell'=1}^d \left|s(k-\ell)s(k'-\ell')s(k-k')s(\ell-\ell')\right| \leq     \frac{1}{d}\left(\sum_{|k|\leq d}|s(k)|^{4/3}\right)^{3},
          \end{align}
          where the second inequality follows from the same computations as in \cite[page 134-135]{NP}.
          
          Moreover, for any $\sigma, \sigma', \tau, \tau' \in S(p)$, 
          \begin{align}\label{2p-1}
         & \sum_{\mathbf{j}, \mathbf{j}' \in \Delta_{p}^{\uparrow} }   
         \left| \left(\prod_{m=1}^{q}r(j_{\sigma(m)} -j'_{\tau(m)})r(j_{\sigma'(m)} -j'_{\tau'(m)})\right)
                \left( \prod_{m=q+1}^{p}r(j_{\sigma(m)} - j_{\sigma'(m)})r(j'_{\tau(m)} - j'_{\tau'(m)})\right)\right|\nonumber\\
                & \quad \leq  \sum_{\mathbf{j}, \mathbf{j}' \in \Delta_{p}^{\uparrow} } \left| r(j_{\sigma(1)}- j'_{\tau(1)})\right|
                \leq n^{2p-2}\sum_{k, \ell=1}^{n}|r(k-\ell)| \leq n^{2p-1}\|r\|_{\ell^1(\Z)}.
          \end{align}
          Therefore,  we combine \eqref{4/3}, \eqref{2p-1} and \eqref{square} to obtain \eqref{2p-1d}.
\end{proof}

We are now at the position to prove Theorem \ref{theorem:tensor}

\begin{proof}[Proof of Theorem \ref{theorem:tensor}]          
           Recall the random vectors $\widetilde{Y}_{n, d}$ and $\mathbf{Z}^{\uparrow}$ defined in \eqref{Yup}. By Lemma \ref{invertible-tensor}, the 
           covariance matrix $\widetilde{C}$ of $\mathbf{Z}^{\uparrow}$ is invertible. Hence,
           according to Proposition \ref{malliavin-stein} and the identity \eqref{expectation}, we have 
           \begin{align}\label{gauss1}
          d_{\text{Wass}}\left(\tilde{\mathcal{Y}}_{n, d}^{\uparrow}\,, \mathbf{Z}^{\uparrow}\right) &\leq \|\widetilde{C}^{-1}\|_{\text{op}}
          \|\widetilde{C}\|^{1/2}_{\text{op}}
          \left(  \sum_{\mathbf{j}, \mathbf{j}' \in \Delta_{p}^{\uparrow} }
          \Var\left(p^{-1}\left\langle DI_p(f^{(d)}_{\mathbf{j}}) \,, DI_p(f^{(d)}_{\mathbf{j}'})\right\rangle_{\eta}\right)
          \right)^{1/2}\nonumber\\
          & = \|\widetilde{C}^{-1}\|_{\text{op}}\|\widetilde{C}\|^{1/2}_{\text{op}}
          \left(  \sum_{\mathbf{j}, \mathbf{j}' \in \Delta_{p}^{\uparrow} }
          O\left(\sum_{q=1}^{p-1}\|f_{\mathbf{j}}^{(d)} \otimes_{q}
            f_{\mathbf{j}'}^{(d)}\|^2_{\eta^{\otimes 2p-2r}} \right)
          \right)^{1/2}\nonumber\\
          & = \|\widetilde{C}^{-1}\|_{\text{op}}\|\widetilde{C}\|^{1/2}_{\text{op}}
          \|r\|^{1/2}_{\ell^1(\Z)}
          \left(\sum_{|k|\leq d}|s(k)|^{4/3}\right)^{3/2}
           O\left(\sqrt{\frac{n^{2p-1}}{d}}\right),
         \end{align}
        where the first equality holds by  \eqref{var} and the second follows from Lemma \ref{squarenorm}. 
        Moreover, we apply Lemma \ref{invertible-tensor} to see that there exists a constant $c_p>0$ such that 
                   \begin{align}\label{gauss2}
          d_{\text{Wass}}\left(\tilde{\mathcal{Y}}_{n, d}^{\uparrow}\,, \mathbf{Z}^{\uparrow}\right) &\leq  c_p
                           \sqrt{\frac{d} {\left|\sum_{k,\ell =1}^d s(k-\ell)^p\right|}}
                           \frac{\left(1 -   \left(\|r\|_{\ell^1(\Z)} -1\right)
               \left(p!  \|r\|^{p-1}_{\ell^1(\Z)} + (p!-1)/2\right)\right)^{1/2}}{ 1 -   \left(\|r\|_{\ell^1(\Z)} -1\right)
              \left(p!  \|r\|^{p-1}_{\ell^1(\Z)} + (p!-1)/2\right)}
             \nonumber\\
               & \qquad \times  \|r\|^{1/2}_{\ell^1(\Z)}
                     \sqrt{\frac{n^{2p-1}}{d} \times \left(\sum_{|k|\leq d}|s(k)|^{4/3}\right)^{3}}.
           \end{align}
           Applying the argument of Step 2 on page 22 of \cite{NourdinZheng2018} (similar to \eqref{full-half}), 
            we conclude the proof of \eqref{gauss}.
\end{proof}

\vskip1cm
\begin{small}
\noindent\textbf{Ivan Nourdin} and \textbf{Fei Pu.}
University of Luxembourg,
Department of Mathematics,
Maison du Nombre,
6 avenue de la Fonte,
L-4364 Esch-sur-Alzette,
Grand Duchy of Luxembourg \\
Emails: \texttt{ivan.nourdin@uni.lu} and \texttt{fei.pu@uni.lu}\\
\end{small}

\end{document}